\documentclass[3p, a4paper, 11pt]{elsarticle}

\usepackage{xcolor}
\usepackage{tikz}
\usepackage{layouts}
\biboptions{numbers,sort&compress}
\usepackage{algorithm}
\usepackage{algpseudocode}
\usepackage{gensymb}
\usepackage[unicode=true,colorlinks=true,pdfborder={0 0 0}]{hyperref}
\hypersetup{
    bookmarks=true,         
    bookmarksopen=true,     
    pdftoolbar=true,        
    pdfmenubar=true,        
    pdffitwindow=true,      
    pdfstartview={FitH},    
    pdfnewwindow=true,      
}
\usepackage{graphicx}
\usepackage{caption}
\usepackage{subcaption}
\usepackage{wrapfig}
\usepackage{multirow,makecell}
\usepackage{longtable}
\usepackage{booktabs}
\usepackage{tabularx}
\usepackage{setspace}
\usepackage{float}
\usepackage{enumerate}
\usepackage{enumitem}

\usepackage[amsthm,thmmarks]{ntheorem}
\usepackage{siunitx}
\usepackage{lineno}

\usepackage{math-symbols}
\usepackage{cases}

\newcommand\etal{{et al.}\ }
\newcommand\ie{{i.e.}}

\graphicspath{{./images/}}

\theoremstyle{definition}

\theoremstyle{definition}
\newtheorem{example}{Example}
\newtheorem{remark}{Remark}
\theoremstyle{definition}
\newtheorem{definition}{Definition}

\theoremstyle{definition}

\journal{XXX} 

\begin{document}
	\begin{frontmatter}
		\title{Balanced truncation for bilinear time-delay systems with approximate 
		Gramians}
		\author[home]{Xiaolong Wang\corref{cor1}}
		\ead{xlwang@nwpu.edu.cn}
		\author[home]{Xin Su}
		\author[home]{Xuerong Yang}
		\cortext[cor1]{Corresponding author}
		\address[home]{School of Mathematics and Statistics,
			Northwestern Polytechnical University, Xi'an 710072, China}
		\begin{abstract}
			A balanced truncation (BT) procedure is presented to achieve a 
			structure-preserving model order reduction (MOR) for bilinear 
			time-delay systems in this 
			paper. We attempt to define Gramians for bilinear time-delay systems 
			based on Volterra series theory. The controllability and observability 
			Gramians in the frequency domain are given explicitly, which align with 
			the ones of bilinear systems and time-delay systems. 
			Based on the derived Gramians, a BT procedure is 
			provided to produce reduced order models with the same structure. We 
			also propose a numerical quadrature rule based on the truncated 
			Laguerrre expansion to give an accurate approximation to Gramians. The 
			resulting low-rank structure of approximate Gramians benefits a lot the 
			efficient execution of the whole MOR procedure. Two numerical examples 
			are simulated to showcase the efficiency of our approach.		
		\end{abstract}
	\end{frontmatter}

	\section{Introduction}\label{sec:introduction}
	Model order reduction (MOR) aims to derive a low-dimensional approximation to 
large-scale systems while preserve the basic dynamical behaviors of original 
systems, thereby reducing dramatically the cost of simulation. A variety of MOR 
methods have been developed in the past decades, such as balanced truncation 
(BT), Krylov subspace-based methods, proper orthogonal decomposition, \etal 
\citep{1102568,pillage2002asymptotic,jiang2011time}.
The fundamental theory on these methods have been well-established for linear 
systems and widely used in engineering field 
\citep{simoncini2007new,benner2008numerical,bai2005dimension,
	salimbahrami2006order,fortuna2012model,lassila2014model,bhattacharya2021model}.
The extension to nonlinear systems has also been considered, and reduced order 
models (ROMs) maintain good accuracy and specific structure under certain 
conditions.

We consider the problem of MOR for bilinear time-delay (BTD)
systems, which are described by the following differential equations
\begin{equation}
	\left\{\begin{aligned}
		\dot{x}(t) & =A_0 x(t)+A_1 x(t-\tau)+\sum_{k=1}^m N_k x(t) u_k(t)+B 
		u(t), \\
		y(t) & =C x(t), t>0, \\
		x(t) & =\phi(t), t \in[-\tau, 0], 
	\end{aligned}\right.
	\label{eq1}
\end{equation}
where $x(t)\in\mathbb{R}^n$ is the state variable, $u(t)=[u_1(t), \cdots, 
u_m(t)]^{\top}\in\mathbb{R}^m$ is the 
input, and $y(t)\in\mathbb{R}^l$ is the output. Further, $\tau$ is the time 
delay, and $A_0, A_1, N_k
\in\mathbb{R}^{n\times n}$. We assume the zero initial function $\phi(t)=0$ in 
this paper. We aim to 
construct 
structure-preserving ROMs of order $r$, such that the 
dynamical behavior of \eqref{eq1} is replicated for admissible inputs 
\begin{equation}
	\left\{\begin{aligned}
		\dot{\hat{x}}(t) & =\hat{A_0} \hat{x}(t)+\hat{A_1} 
		\hat{x}(t-\tau)+\sum_{k=1}^m \hat{N_k} \hat{x}(t) \hat{u_k}(t)+\hat{B} 
		\hat{u}(t), \\
		\hat{y}(t) 	& = \hat{C} \hat{x}(t), t>0, \\
		\hat{x}(t) & =\hat{\phi}(t), t \in[-\tau, 0],
	\end{aligned}\right.
	\label{eq2}
\end{equation}
where $\hat{x} \in \mathbb{R}^r, \hat{A_0}, \hat{A_1}, \hat{N_k} \in 
\mathbb{R}^{r \times r} (k = 1,2,\cdots, m), \hat{B} \in \mathbb{R}^{r \times 
	m}$ 
and $\hat{C} \in \mathbb{R}^{l \times r}$ for $r \ll n.$ BTD systems take on the 
properties of bilinear systems and time-delay systems simultaneously.

Time delays are ubiquitous in practical, arising frequently in problems such 
as signal 
propagation in communication networks, decision feedback in economic systems, 
and neural signal transmission in biomedicine \citep{erneux2009applied, 
	kuang1993delay, fridman2014introduction, smith2011introduction}.
Conventional methods that do not account for the effect coming from time 
delays often perform poorly in MOR of large-scale time-delay systems. 
Consequently, the study of MOR for time-delay systems carries 
significant theoretical and practical value. Several approaches on MOR of 
time-delay systems	are available in the literature. 
Michiels \etal transformed the time-delay system into 
an equivalent infinite-dimensional linear system \citep{michiels2011krylov}. By 
combining spectral 
discretization with the Padé approximation, the approach achieves moment 
matching at the infinity and the origin points, but cannot preserve the time 
delay 
structure during MOR.
Kharitonov and Zhabko introduced a couple of delay Lyapunov matrices for the 
sake of stability analysis, which also established the
theoretical foundation of Gramians \citep{kharitonov2003lyapunov, 
	kharitonov2006lyapunov} for time-delay systems.
Later, a Gramian-based MOR technique was developed for linear time-delay systems 
based on positional balancing, which preserves the delay structure of systems 
\citep{jarlebring2013model}.
Another BT method is also provided for asymptotically stable time-delay 
systems in \citep{8264217}.
Lordejani \etal applied the BT procedure to a special case where the output 
depends on delays as well, and the structure and stability of original systems  
can be guaranteed during MOR \citep{lordejani2020model}. Wang and 
Jiang derived a frequency-domain approximation to
controllability and observability Gramians, thereby enabling an efficient 
execution of BT for time-delay systems \citep{wang2011balanced}.
Xiao and Jiang investigated Krylov subspace-based  MOR for second-order 
time-delay systems, and the two types of structure-preserving MOR schemes are 
capable of matching a prescribed number of moments or Laguerre 
coefficients \citep{xiao2016multi}.

Bilinear systems are a class of special nonlinear systems whose dynamics are 
governed by differential equations containing bilinear terms.
A amount of  nonlinear dynamical systems can be well approximated via bilinear 
systems through Carleman linearization.
Consequently, the development of MOR for bilinear systems also serves as an 
important approach to achieve the reduction of nonlinear dynamical models.
Zhang and Lam investigated the necessary conditions of $\mathcal{H}_2$-optimal 
MOR for continuous-time bilinear 
systems \citep{zhang2002h2}.
Wang and Jiang constructed ROMs based on Laguerre series expansion and 
investigated 
the relationship to the typical moment-matching methods \citep{wang2012model}.
Antoulas \etal generalized the Loewner framework from linear systems to bilinear 
systems, which allows for the direct construction of ROMs 
based on the 
input-output 
sample data, without the need to access the full order models  
\citep{antoulas2016model}.
Nasiri Soloklo \etal proposed an $\mathcal{H}_2$-optimal MOR method based 
on linear matrix inequalities, and the stability of systems is ensured by 
imposing the constraint on coefficient matrices of ROMs directly
\citep{nasiri2023h2}.

This paper focuses on the BT approach of BTD systems. As this type of 
systems couples both the nonlinearity from bilinear terms and the delay effect 
from time delay structure, MOR methods of BTD systems remain rather limited 
in the literature.
Gosea \etal adopted the Loewner framework to match the moments of multivariate 
transfer functions for BTD systems \citep{gosea2019model}.
A Krylov subspace-based MOR method was proposed for BTD systems based on
Taylor and Laguerre series expansion of multivariable transfer functions, and a 
desired 
number of moments can be preserved by ROMs \citep{cheng2024high}.
Currently, there is no Gramian-based methods available for MOR of BTD systems. 
We attempt to establish the concept of Gramians for BTD systems. Specifically, 
Gramians are first defined in the time domain based on the Volterra series 
theory, and then the explicit expression of Gramians is provided in the 
frequency domain via Laplace transformation. As a consequence, a BT procedure 
is presented to produce structure-preserving ROMs relying on the derived 
Gramians. Furthermore, we also propose a numerical quadrature
rule with the aid of truncated Laguerrre expansions to obtain a low-rank  
approximation to Gramians, which leads to an efficient execution of our 
approach.

The remainder of the paper is structured as follows. 
The preliminary on linear time-delay systems and bilinear systems is reviewed 
in Section 2. The definition of Gramians for BTD systems is presented first 
based on the Volterra series theory in Section 3, followed by a BT procedure 
used to 
produce structure-preserving ROMs. The low-rank approximation to Gramians is 
also provided in Section 3 to achieve an efficient execution of our approach.  
Section 4 presents the simulation results based on two numerical examples. A 
summary 
and possible topics of future works are given in Section 5.

\section{Preliminary}\label{sec:sec-2}
We review some preliminaries on bilinear systems and time-delay 
systems in this section, which pave the way for defining the basic concepts of 
BTD systems.

Linear time-delay systems with a single delay are describes as 
\begin{equation}
	\left\{\begin{aligned}
	\dot{x}(t) & =A_0 x(t)+A_1 x(t-\tau)+B u(t), \\
	y(t) & =C x(t), t>0, \\
	x(t) & =\phi(t), t \in[-\tau, 0], 
	\end{aligned}\right.
	\label{eq3}
	\end{equation}
where $x(t) \in \mathbb{R}^n$, $u(t) \in \mathbb{R}^m$, and $y(t) \in 
\mathbb{R}^l$ represent the system state, input, and output, respectively.
$\tau$ is the delay, and 
$\phi(t)$ is the given initial function.

\begin{definition}\label{def1}
	For a given linear time-delay system \eqref{eq3}, the fundamental solution 
	$K(t) \in \mathbb{R}^{n 
	\times n}$ of \eqref{eq3} is the solution to the matrix delay-differential 
	equation
	\begin{equation*}
		\left\{\begin{array}{l}
		\dot{K}(t)=A_0 K(t)+A_1 K(t-\tau), \\
		K(0)=I, K(\theta)=0, \theta<0.
		\end{array}\right.
	\end{equation*}
\end{definition}
For an exponentially stable linear time-delay system \eqref{eq3}, $K(t)$ decays 
exponentially, and the associated time-delay Lyapunov matrix can be 
defined \citep{jarlebring2013model}.

\begin{definition}\label{def2}
	For a given constant $\theta$, the delay Lyapunov matrices 
	related to controllability and observability
	are defined, respectively, as		
	\begin{equation*}
			\begin{aligned}
			& U_c(\theta)=\int_0^{\infty} K(t) B B^{\top} K^{\top}(t+\theta) d t, \\
			& U_o(\theta)=\int_0^{\infty} K^{\top}(t) C^{\top} C K(t+\theta) d t .
			\end{aligned}
	\end{equation*}
	One can verify that $U_c: [-\tau,\tau]\rightarrow\mathbb{R}^{n\times n}$ and 
	$U_o: [-\tau,\tau]\rightarrow\mathbb{R}^{n\times n}$ are 
	respectively the unique solutions of the following boundary value problems
	\begin{equation*}
			\begin{aligned}
			& \left\{\begin{array}{l}
			\dot{U}_c(t)=U_c(t) A_0^{\top}+U_c(t-\tau) A_1^{\top}, t \geqslant 0, \\
			U_c(-t)=U_c(t)^{\top}, \\
			-B B^{\top}=U_c(0) A_0^{\top}+A_0^{\top} U_c(0)+U_c(-\tau) 
			A_1^{\top}+A_1 U_c(\tau) .
			\end{array}\right. \\
			& \left\{\begin{array}{l}
			\dot{U}_o(t)=U_o(t) A_0+U_o(t-\tau) A_1, t \geqslant 0, \\
			U_o(-t)=U_o(t)^{\top}, \\
			-C^{\top} C=U_o(0) A_0+A_0^{\top} U_o(0)+U_o(-\tau) A_1+A_1^{\top} 
			U_o(\tau) .
			\end{array}\right.
			\end{aligned}
	\end{equation*}
	In particular, $U_c(0)$ and $U_o(0)$ are defined as the controllability and 
	observability Gramians of\eqref{eq3}, and they have the following explicit 
	expression
	\begin{equation*}
			\begin{aligned}
			& U_c(0)=\int_0^{\infty} K(t) B B^{\top} K^{\top}(t) d t, \\
			& U_o(0)=\int_0^{\infty} K^{\top}(t) C^{\top} C K(t) d t.
			\end{aligned}
	\end{equation*}
	
	\end{definition}

	Bilinear systems are a special type of nonlinear systems, which can be 
	described as follows 
		\begin{equation}
			\left\{\begin{aligned}
			\dot{x}(t) & =A_0 x(t)+\sum_{k=1}^m N_k x(t) u_k(t)+B u(t), \\
			y(t) & =C x(t), t>0, \\
			x(0) & =0,
			\end{aligned}\right.
			\label{eq4}
		\end{equation}
	where the nonlinearity comes from the interplay in the bilinear terms, and 
	the coefficient matrix $N_k \in \mathbb{R}^{n \times n} (k = 1, 2, \cdots, 
	m)$. The balanced truncation procedure has been studied extensively for 
	\eqref{eq4}, and Gramians are well defined for stable bilinear systems 
	\cite{Benner2011add}. 
	
	\begin{definition}\label{def3}
		The controllability and observability Gramians of \eqref{eq4} are  
		defined, respectively, as		  
		\begin{equation*}
		\begin{aligned}
  \bar{P} & =\sum_{j=1}^{+\infty} \int_{0}^{+\infty} \cdots 
  \int_{0}^{+\infty} \bar{P}_j\left(t_1, \cdots,  t_j\right) 
  \bar{P}_j^\top\left(t_1, \cdots,  t_j\right) d t_1 \cdots d t_j, \\
  \bar{Q} & =\sum_{j=1}^{+\infty} \int_{0}^{+\infty} \cdots 
  \int_{0}^{+\infty} \bar{Q}_j^\top\left(t_1, \cdots,  t_j\right) 
  \bar{Q}_j\left(t_1, \cdots,  t_j\right) d t_1 \cdots d t_j,
		\end{aligned}
		\end{equation*}
		where $\bar{P}_j\left(t_1, \cdots,  t_j\right)$ and $\bar{Q_j}\left(t_1, 
		\cdots,  t_j\right)$ for $j = 1,2, \cdots $ have the explicit form 
		\begin{equation*}
  \begin{aligned}
  & \bar{P}_1\left(t_1\right) = e^{A t_1} B, \quad
  \bar{P}_j\left(t_1, \cdots,  t_j\right)=e^{A t_j}\left[N_1 P_{j-1}, 
  \cdots, N_m P_{j-1}\right], \\
  & \bar{Q}_1\left(t_1\right)=Ce^{A t_1}, \quad 
  \bar{Q}_j\left(t_1, \cdots,  t_j\right)=\left[\begin{array}{c}
  Q_{j-1} N_1 \\
  \vdots \\
  Q_{j-1} N_m
  \end{array}\right]e^{A t_j}. 
  \end{aligned}
  \end{equation*}
 \end{definition}
Gramains of bilinear systems satisfy a couple of generalized Lyapunov equtions, 
and can be calculated via the low-rank approximation. Gramians play an important 
role in the standard BT procedure. We will define Gramians for BTD systems 
in the next section following the same lines. 


\section{BT procedure for BTD systems}\label{sec:sec-3}
In this section, we first define Gramians for BTD systems based on the Volterra 
representation, and then present a BT procedure to produce structure-preserving 
ROMs. The approximate Gramians based on Laguerre expansion are also introduced 
to enable an efficient execution of our approach. 

\subsection{Gramians of BTD systems}

A basic tool for the analysis of BT methods is Volterra series expansion of 
dynamical systems. We now define Gramians of BTS systems based on the theory of 
Volterra series expansion \cite{rugh1981nonlinear}. Assume that the state $x(t)$ 
of \eqref{eq1} possesses the following asymptotic expansion with respect to the  
parametric input $\alpha u(t)$
\begin{equation}
	x(t)=\sum_{i=1}^{\infty} \alpha^i x_i(t),
	\label{eq5}
\end{equation}
where $\alpha\in\mathbb{R}$ is a parameter, and $x_i(t)\in\mathbb{R}^n$. Taking 
it 
into the state equation of 
\eqref{eq1} along with the input $\alpha u(t)$ leads to
\begin{equation*}
	\begin{aligned}
	\sum_{i=1}^{\infty} \alpha^i \dot{x}_i(t)= & A_0\left(\sum_{i=1}^{\infty} 
	\alpha^i x_i(t)\right) + A_1 \left(\sum_{i=1}^{\infty} \alpha^i 
	x_i(t-\tau)\right)\\
	& +\sum_{k=1}^m \sum_{i=1}^{\infty}N_k\left(\alpha^{i+1} x_i(t)\right) 
	u_k(t)+\alpha B u(t).
	\end{aligned}
	\label{eq6}
\end{equation*}
As $\alpha$ is an arbitrary parameter, equating the coefficients of $\alpha^i$ 
in both sides of the above equation results in a series of differential 
equations 
\begin{equation*}
	\begin{aligned}
	& \dot{x}_1(t)=A_0 x_1(t) + A_1 x_1(t-\tau) + B u(t), \\
	& \dot{x}_i(t)=A_0 x_i(t) + A_1 x_i(t-\tau) + \sum_{k=1}^m N_k x_{i-1}(t) 
	u_k(t), \quad i \geq 2 .
	\end{aligned}
\end{equation*}
Note that the equations of $x_i(t)$ are standard linear time-delay differential 
equations coupled in a sequential order. With the zero initial conditions for 
the above equations, $x_1(t)$ 
can be expressed via the 
fundamental solution in Definition 1 as follows 
\begin{equation*}
	x_1(t)=\int_0^t K(t_1) B u\left(t-t_1\right) d t_1.
\end{equation*}
Let $P_1\left(t_1\right) = K(t_1) B$ be the first kernel function. $x_2(t)$ is 
described 
in the same 
manner 
\begin{equation*}
	\begin{aligned}
		x_2(t)	&=\sum_{k=1}^m \int_0^{t} \int_0^{t-t_2} K(t_2) N_k K(t_1) B 
		u(t-t_2-t_1) u_k\left(t-t_2\right) d t_1 d t_2.\\
				&= \int_0^t \int_0^{t-t_2} K(t_2)\left[N_1 P_1(t_1), \cdots, N_m 
				P_1(t_1)\right] (u(t-t_2) \otimes u\left(t-t_1-t_2\right)) d t_1 
				d t_2,
	\end{aligned}
 	\end{equation*}
where $\otimes$ represents the Kronecker product. In general, with the 
$(i-1)$-th kernel function $$P_{i-1}\left(t_1, \cdots,  
t_{i-1}\right)=K(t_{i-1})\left[N_1 
P_{i-2}\left(t_1, \cdots,  t_{i-2}\right), 
\cdots, N_mP_{i-2}\left(t_1, \cdots,  t_{i-2}\right)\right], $$ 
we can obtain a unified expression for $x_i(t)$ 
\begin{equation*}
   	\begin{aligned}
   		x_i(t)=	&\int_0^t \int_0^{t-t_i} \cdots \int_0^{t-t_i-\cdots-t_2} 
   		K(t_i)\left[N_1 P_{i-1}\left(t_1, \cdots,  t_{i-1}\right), \cdots, N_m 
   		P_{i-1}\left(t_1, \cdots,  t_{i-1}\right)\right] \\
   				&\times \left(u(t-t_i) \otimes \cdots \otimes 
   				u\left(t-t_1-\cdots-t_i\right)\right)d t_1 \cdots d t_i.
   	\end{aligned}
\end{equation*}
Let $\alpha = 1$. It follows from \eqref{eq5} that $x(t) = \sum_{i=1}^{\infty} 
x_i(t)$, and the explicitly asymptotic expansion is as follows 
\begin{equation*}
	\begin{aligned}
		x(t)= & \int_0^t  K(t_1) B u\left(t-t_1\right) d t_1 \\
		& + \int_0^t \int_0^{t-t_2}  K(t_2)\left[N_1 P_1(t_1), \cdots, N_m 
		P_1(t_1)\right] (u(t-t_2) \otimes u\left(t-t_1-t_2\right)) d t_1 d t_2+ 
		\cdots\\
		& + \int_0^t \int_0^{t-t_i} \cdots \int_0^{t-t_i-\cdots-t_2} \left[N_1 
		P_{i-1}\left(t_1, \cdots,  t_{i-1}\right), \cdots, N_m P_{i-1}\left(t_1, 
		\cdots,  t_{i-1}\right)\right] \\
		& \times \left(u(t-t_i) \otimes \cdots \otimes 
		u\left(t-t_1-\cdots-t_i\right)\right)d t_1 \cdots d t_i+ \cdots.
	\end{aligned}
\end{equation*}

The kernel functions of the asymptotic expansion of \eqref{eq1} paves the way 
for the controllability Gramian of BTD systems, while the 
observability Gramian can be defined similarly via the dual property, as shown 
in Definition \ref{def4}. 

\begin{definition}\label{def4}
	Given BTD systems \eqref{eq1}, the controllability and observability 
	Gramians are defined, respectively, as		
		\begin{equation}\label{eq7}
			\begin{aligned}
				P & =\sum_{i=1}^{+\infty} \int_{0}^{+\infty} \cdots 
				\int_{0}^{+\infty} P_{i}\left(t_1, \cdots,  
				t_{i}\right) P_{i}\left(t_1, \cdots,  
				t_{i}\right)^\top d t_1 \cdots d t_i, \\
				Q & =\sum_{i=1}^{+\infty} \int_{0}^{+\infty} \cdots 
				\int_{0}^{+\infty} Q_i\left(t_1, \cdots,  t_i\right)^\top 
				Q_i\left(t_1, \cdots,  t_i\right) d t_1 \cdots d t_i,
			\end{aligned}
		\end{equation}
	where the kernel functions are given by
		\begin{equation*}
			\begin{aligned}
				& P_1\left(t_1\right) = K(t_1) B, \\
				& P_i\left(t_1, \cdots,  t_i\right)=K(t_i)\left[N_1 
				P_{i-1}\left(t_1, \cdots,  
				t_{i-1}\right), 
				\cdots, N_m P_{i-1}\left(t_1, \cdots,  
				t_{i-1}\right)\right], \\
				& Q_1\left(t_1\right)=CK(t_1), \\
				& Q_i\left(t_1, \cdots,  t_i\right)=\left[\begin{array}{c}
				Q_{i-1}\left(t_1, \cdots,  
				t_{i-1}\right) N_1 \\
				\vdots \\
				Q_{i-1}\left(t_1, \cdots,  
				t_{i-1}\right) N_m
				\end{array}\right]K(t_i), i = 2,3, \cdots .
			\end{aligned}
		\end{equation*}
	in which $ K(t_i)$ is the fundamental solution in Definition \ref{def1}.
\end{definition}

In order to calculate Gramians efficiently, we perform a zero continuation 
for $P_i(t_1, \cdots, t_i)$ and $Q_i(t_1, \cdots, t_i)$ in the interval 
$t_i\in[-\infty, +\infty]$ and consider the manipulation in the frequency 
domain. Based on Parseval theorem, the Laplace transform of $P_i(t_1, 
\cdots, t_i)$ and $Q_i(t_1, \cdots, t_i)$ leads to 		
\begin{equation}\label{eq8}
	\begin{aligned}
		& P=\sum_{i=1}^{\infty} \frac{1}{(2 \pi)^i} \int_{-\infty}^{+\infty} 
		\cdots \int_{-\infty}^{+\infty} P_i\left(\mathrm i w_1, \cdots, 
		\mathrm i w_i\right) 
		P_i\left(\mathrm i w_1, \cdots, \mathrm i w_i\right)^{\mathrm H} d 
		w_1 \cdots d 
		w_i,\\
		& Q=\sum_{i=1}^{\infty} \frac{1}{(2 \pi)^i} \int_{-\infty}^{+\infty} 
		\cdots \int_{-\infty}^{+\infty} Q_i\left(\mathrm i w_1, \cdots, 
		\mathrm i 
		w_i\right)^{\mathrm H} Q_i\left(\mathrm i w_1, \cdots, \mathrm i 
		w_i\right) d w_1 
		\cdots d 
		w_i,
	\end{aligned}
\end{equation}
where $\mathrm i = \sqrt{-1}$, and the integrands are given by 
\begin{equation*}
	\begin{aligned}
		& P_1\left(\mathrm i w_1\right)=\left(\mathrm i w_1 
		I-A_0-e^{-\mathrm i w_1 \tau} A_1\right)^{-1} B,\\
		& P_i\left(\mathrm i w_1, \cdots, \mathrm i w_i\right)=\left(\mathrm 
		i w_i I-A_0-e^{-\mathrm i w_i \tau} 
		A_1\right)^{-1}\left[\begin{array}{lll}
		N_1 P_{i-1}, & \cdots, & N_m P_{i-1}
		\end{array}\right], \quad i=2,3, \cdots ,  \\
		& Q_1\left(\mathrm i w_1\right)=C\left(\mathrm i w_1 
		I-A_0-e^{-\mathrm i w_1 \tau} A_1\right)^{-1}, \\
		& Q_i\left(\mathrm i w_1, \cdots, \mathrm i 
		w_i\right)=\left[\begin{array}{c}
		Q_{i-1} N_1 \\
		Q_{i-1} N_2 \\
		\cdots \\
		Q_{i-1} N_m
		\end{array}\right]\left(\mathrm i w_i I-A_0-e^{-\mathrm i w_i \tau} 
		A_1\right)^{-1}, \quad i=2,3, \cdots .
	\end{aligned}
\end{equation*}

\begin{remark}
Obviously, when the coefficient matrix $A_1 = 0$, BTD systems boil down to the 
standard bilinear systems, and Gramians defined in \eqref{eq8} coincide with ones 
given in \cite{Benner2011add}. In addition, if $N_k=0$ for $k=1, \cdots, m$, the 
Gramians of BTD systems degenerate into ones of linear 
time-delay systems in \cite{jarlebring2013model}. 
\end{remark}

\subsection{Approximate Gramians based on Laguerre expansion}
Although Gramians of BTD systems are given explicitly in \eqref{eq7} and 
\eqref{eq8}, the involved quadrature and infinite series are much more complex and 
cannot be calculated efficiently in practice. In this subsection, we provide a 
scheme to approximate Gramians based on a Laguerre expansion strategy.

Laguerre polynomials in the Rodrigues representation form are given by
	\begin{equation*}
		L_i(t)=\frac{e^t}{i!} \frac{d^i}{d t^i}\left(e^{-t} t^i\right), \quad i=0,1, 
		\cdots,
	\end{equation*}
and the associated Laguerre function is	adopted to enable an expansion in our 
settings 
	\begin{equation*}
		\phi_i^\alpha(t)=\sqrt{2 \alpha} e^{-\alpha t} L_i(2 \alpha t), \quad i=0,1, 
		\cdots,
	\end{equation*}
where $\alpha > 0$ is referred as the time scale factor \citep{wang2012model}. It is 
well known that Laguerre functions form a uniformly bounded orthonormal basis for 
the Hilbert space $L_2(\mathbb R_+)$. 
Laplace transformation of Laguerre functions leads to the following frequency domain 
expression
	\begin{equation*}
		\Phi_i^\alpha(s)=\frac{\sqrt{2 
		\alpha}}{s+\alpha}\left(\frac{s-\alpha}{s+\alpha}\right)^i, \quad i=0,1, 
		\cdots,
	\end{equation*}
which satisfies the orthogonal property	
	\begin{equation}\label{orthogonal_p}
		\frac{1}{2 \pi} \int_{-\infty}^{+\infty} \Phi_i^\alpha(\mathrm i 
		\omega)\left(\Phi_j^\alpha(\mathrm i \omega)\right)^{\mathrm{H}} d 
		\omega=\delta_{ij},
	\end{equation}
where $\delta_{k q}$ is the Kronecker function.

We consider the Laguerre expansion of $\left(s I - A_0 - e^{-s \tau} 
A_1\right)^{-1}$, which is the focus of \eqref{eq7} and \eqref{eq8}. 
We assume that 
\begin{equation}\label{eq10}
	\left(s I-A_0-e^{-s \tau} A_1\right)^{-1}=\sum_{i=0}^{\infty} \tilde{f}_i 
	\Phi_i^\alpha(s), 
\end{equation}
where $\tilde{f}_i\in\mathbb{R}^{n\times n}$ are the Laguerre coefficients. 
By introduce the transformation $u = (s - \alpha)/(s + \alpha)$, there holds
$$
\Phi_i^\alpha(s)=\frac{1-u}{\sqrt{2 \alpha}} u^i.
$$
The expansion \eqref{eq10} is reformulated as  
\begin{equation*}
	\left(s I-A_0-e^{-s \tau} A_1\right)^{-1}=\sum_{i=0}^{\infty} \tilde{f}_i 
	\Phi_i^\alpha(s)=\sum_{i=0}^{\infty} \tilde{f}_i \frac{1-u}{\sqrt{2 \alpha}} 
	u^i.
\end{equation*}
Since $\left(s I-A_0-e^{-s \tau} A_1\right)\left(s I-A_0-e^{-s \tau} 
A_1\right)^{-1}=I$, a simple manipulation leads to 
	\begin{equation}\label{eq12}
		\left(\alpha(1+u) I-(1-u) A_0-(1-u) e^{-\frac{\alpha(1+u)}{1-u} \tau} 
		A_1\right) \sum_{i=0}^{\infty} \tilde{f}_i u^i=\sqrt{2 \alpha}I.
	\end{equation}
Because the generating function of Laguerre polynomials is $$\frac{e^{-\frac{x 
			t}{1-t}}}{1-t}=\sum_{i=0}^{\infty} L_i(x) t^i,$$
the exponential function in \eqref{eq12} can be expanded as 
\begin{equation*}
		e^{-\frac{\alpha(1+u)}{1-u} \tau}=e^{-\alpha \tau} \sum_{i=0}^{\infty} L_i(2 
		\alpha \tau)\left(u^i-u^{i+1}\right).
\end{equation*}
Substituting it into \eqref{eq12}, we get
\begin{equation*}
		\left(\alpha(1+u) I-(1-u) A_0-(1-u) e^{-\alpha \tau} \sum_{i=0}^{\infty} 
		L_i(2 \alpha \tau)\left(u^i-u^{i+1}\right) A_1\right) \sum_{i=0}^{\infty} 
		\tilde{f}_i u^i=\sqrt{2 \alpha} I.
\end{equation*}
With the notations 
	\begin{equation}\label{Delata_Coeff}
	\begin{aligned}
		& \triangle_0=\alpha I-A_0-L_0(2 \alpha \tau) e^{-\alpha \tau} A_1, \\
		& \triangle_1=\alpha I+A_0-\left(-2 L_0(2 \alpha \tau)+L_1(2 \alpha 
		\tau)\right) e^{-\alpha \tau} A_1, \\
		& \triangle_i=-\left(L_{i-2}(2 \alpha \tau)-2 L_{i-1}(2 \alpha 
		\tau)+L_i(2 \alpha \tau)\right) e^{-\alpha \tau} A_1, \quad i \geqslant 2, 
	\end{aligned}
\end{equation}
the above equation is formulated as a compact form 
\begin{equation*}\label{eq13}
	\sum_{i=0}^{\infty} \sum_{j=0}^{\infty} \Delta_i \tilde{f}_j u^{i+j}=\sqrt{2 
		\alpha} I. 
\end{equation*}
By equating the coefficients of $u^i$ for $i=0, 1, \cdots, q-1, $ we obtain the 
following linear systems about Laguerre coefficients 
	\begin{equation*}
		\left[\begin{array}{cccc}
		\Delta_0 & & & \\
		\Delta_1 & \Delta_0 & & \\
		\vdots & \ddots & \ddots & \\
		\Delta_{q-1} & \cdots & \Delta_{1} & \Delta_0
		\end{array}\right]\left[\begin{array}{c}
		\tilde{f}_0 \\
		\tilde{f}_1 \\
		\vdots \\
		\tilde{f}_{q-1}
		\end{array}\right]=\left[\begin{array}{c}
		\sqrt{2 \alpha} I \\
		0 \\
		\vdots \\
		0
		\end{array}\right].
	\end{equation*}
If the coefficient matrix $\Delta_0$ is invertible, Laguerre coefficients 
$\tilde{f}_i$ are determined uniquely. 
Due to a block lower triangular structure of the coefficient matrix, $\tilde{f}_i$ 
can be obtained recursively 
\begin{equation}\label{coeff_expp}
	\begin{aligned}
		\tilde{f}_0 & =\triangle_0^{-1} \sqrt{2 \alpha} I, \\
		\tilde{f}_1 & =-\triangle_0^{-1}\left(\triangle_1 \tilde{f}_0\right), \\
		\tilde{f}_i & =-\triangle_0^{-1}\left(\triangle_i 
		\tilde{f}_0+\triangle_{i-1} \tilde{f}_1+\cdots+\triangle_1 
		\tilde{f}_{i-1}\right), \quad i=2, \cdots, q-1.
	\end{aligned}
\end{equation}

\begin{remark}
	The invertibility of $\Delta_0$ cannot be ensured naturally. One can adjust the 
	parameter $\alpha$ properly to produce an invertible $\Delta_0$.
	Furthermore, Laguerre coefficients can be 
	computed through the LU decomposition of $\Delta_0$, instead of explicitly 
	calculating the inverse of $\Delta_0$ in practice.
\end{remark}

We now focus on the low-rank approximation to Gramians based on Laguerre expansion.  
For the first term in the right side of \eqref{eq8}, the integrand 
$P_1\left(\mathrm i\omega_1\right)$ has the expansion 
\begin{equation}\label{app_p1}
	P_1\left(\mathrm i\omega_1\right)=\sum_{i=0}^{\infty} \tilde{f}_iB
	\Phi_i^\alpha(\mathrm i\omega_1)=\sum_{i=0}^{\infty} f_{1, i} 
	\Phi_i^\alpha(\mathrm i\omega_1) 
	\approx 
	\sum_{i=0}^{q-1} f_{1, i} \Phi_i^\alpha(\mathrm i\omega_1).
\end{equation}
Note that there holds $f_{1, i}=\tilde f_iB$ for $i=0, 1, \cdots$. 
The above approximation to $P_1\left(\mathrm i\omega_1\right)$ combined with the 
orthogonal 
property \eqref{orthogonal_p} leads to
\begin{equation*}
	\begin{aligned}
		\frac{1}{2 \pi} \int_{-\infty}^{+\infty} P_1\left(\mathrm i\omega_1\right) 
		P_1\left(\mathrm i\omega_1\right)^{\mathrm{H}} d w_1 & 
		\approx \frac{1}{2 \pi} \int_{-\infty}^{+\infty}\left(\sum_{i=0}^{q-1} 
		f_{1, i} \Phi_i\left(\mathrm i w_1\right)\right)\left(\sum_{i=0}^{q-1} f_{1, 
		i} 
		\Phi_i\left(\mathrm i w_1\right)\right)^{\mathrm{H}} d w_1 \\
		& =\sum_{i=0}^{q-1} f_{1, i} f_{1, i}^{\top}=F_1 F_1^{\top}, 
	\end{aligned}
\end{equation*}
where $F_1$ is given by 
\begin{equation}\label{F1}
	F_1 = [f_{1,0}, f_{1,1}, \cdots, f_{1,q-1}], \quad f_{1, i}=\tilde f_iB. 
\end{equation}
Similarly, with the approximation in \eqref{app_p1}, the integrand for the second 
term in \eqref{eq8} boils down to 
\begin{equation*}
	\left(\mathrm i\omega_2 I-A_0-e^{-\mathrm i\omega_2 \tau} 
	A_1\right)^{-1}\left[N_1 F_1, \cdots, N_m 
	F_1\right]=\sum_{i=0}^{\infty} \tilde{f}_i \left[N_1 
	F_1, \cdots, 
	N_m F_1\right]\Phi_q(\mathrm i\omega_2)\approx\sum_{i=0}^{q-1} f_{2, i} 
	\Phi_i(\mathrm i\omega_2), 
\end{equation*}
where there holds $f_{2, i}=\tilde{f}_i \left[N_1 
F_1, \cdots, N_m F_1\right]$ for $i=0, 1, \cdots$. 
The second term of \eqref{eq8} can be approximated as
$$
	\frac{1}{(2 \pi)^2} \int_{-\infty}^{+\infty} \int_{-\infty}^{+\infty} 
	P_2\left(\mathrm i w_2, \mathrm i w_1\right) 
	P_2\left(\mathrm i w_2, \mathrm i w_1\right)^{\mathrm{H}} d w_1 d w_2 \approx 
	\sum_{i=0}^{q-1} f_{2, i} f_{2, 
	i}^{\top}=F_2 F_2^{\top},
$$
where $F_2 = [f_{2,0}, f_{2,1}, \cdots, f_{2,q-1}]$. 
In general, one can show that the $j$-th term of \eqref{eq8} has the 
low-rank approximation 
$$
	\frac{1}{(2 \pi)^j} \int_{-\infty}^{+\infty} \cdots \int_{-\infty}^{+\infty} 
	P_j\left(\mathrm i w_j,\cdots, \mathrm i w_1\right) P_j\left(\mathrm i 
	w_j,\cdots, \mathrm i w_1\right)^{\mathrm{H}} d w_1 \cdots d w_j \approx F_j 
	F_j^{\top},
$$
and the low-rank factor $F_j$ is determined by 
\begin{equation}\label{coeff_pf}
	F_j=[f_{j,0},f_{j,1},\cdots,f_{j,q-1}], \quad f_{j,i}=\tilde{f}_i \left[N_1 
	F_{j-1}, \cdots, N_m F_{j-1}\right]
\end{equation}
for $j=2, 3, \cdots, p$ and $i=1, 2, \cdots, q$. 
Therefore, the low-rank approximation to the controllability Gramian $P$ reads 
\begin{equation*}
		P \approx F_1 F_1^{\top}+F_2 F_2^{\top}+\cdots+F_p F_p^{\top}.
\end{equation*}

As for the observability Gramian $Q$ of BTD systems, we consider the expansion 
\begin{equation*}
	\left(s I-A_0^\top-e^{-s \tau} A_1^\top\right)^{-1}=\sum_{i=0}^{\infty} 
	\tilde{g}_i 
	\Phi_i^\alpha(s). 
\end{equation*}
It follows from \eqref{eq10} that 
\begin{equation}\label{coeff_expQ}
	\begin{aligned}
		\tilde{g}_0 & =\triangle_0^{-\top} \sqrt{2 \alpha} I, \\
		\tilde{g}_1 & =-\triangle_0^{-\top}\left(\triangle_1^\top 
		\tilde{g}_0\right), \\
		\tilde{g}_i & =-\triangle_0^{-\top}\left(\triangle_i^\top 
		\tilde{g}_0+\triangle_{i-1}^\top \tilde{g}_1+\cdots+\triangle_1^\top 
		\tilde{g}_{i-1}\right), \quad i=2, \cdots, q-1.
	\end{aligned}
\end{equation}
With the truncated approximation $Q_1\left(\mathrm i w_1\right)^{\mathrm H}\approx 
\sum_{i=0}^{q-1} g_{1, i} \Phi_i^\alpha(\mathrm i\omega_1)$, the first term of 
\eqref{eq8} has the low-rank approximation 
\begin{equation*}
	\frac{1}{(2 \pi)} \int_{-\infty}^{+\infty}
	Q_1\left(\mathrm i w_1\right)^{\mathrm H} Q_1\left(\mathrm i 
	w_1\right) d w_1 \approx G_1 G_1^{\top},
\end{equation*}
where the factor is determined by 
\begin{equation}\label{G1}
	G_1=[g_{1,0}, g_{1,1}, \cdots, g_{1, q-1}],\quad g_{1,i}=\tilde 
	g_iC^\top. 
\end{equation}
The $j$-th term of the summation in \eqref{eq8} has the approximation 
\begin{equation*}
	\frac{1}{(2 \pi)^j} \int_{-\infty}^{+\infty} \cdots \int_{-\infty}^{+\infty} 
	Q_j\left(\mathrm i w_1, \cdots, \mathrm i 
	w_j\right)^{\mathrm{H}} Q_j\left(\mathrm i w_1, \cdots, \mathrm i 
	w_j\right) d w_1 \cdots d w_j \approx G_j G_j^{\top},
\end{equation*}
where the low-rank factor is given by 
\begin{equation}\label{Gi}
	G_i=[g_{j,0}, g_{j,1}, \cdots, g_{j, q-1}], \quad g_{j,i}=\tilde 
	g_i\left[N_1^\top G_{j-1}, \cdots, N_m^\top G_{j-1}\right] 
\end{equation}
for $j=2, 3, \cdots, p$ and $i=1, 2, \cdots, q$. 
Therefore, the low-rank approximation to the observability Gramian $Q$ is derived 
\begin{equation*}
Q \approx G_1 G_1^{\top}+G_2 G_2^{\top}+\cdots+G_p G_p^{\top}.
\end{equation*}

We summarize the main steps of producing a low-rank approximation to Gramians of BTD 
systems in Algorithm \ref{alg2}.
 
\begin{algorithm}[h]
	\caption{Low-rank approximate Gramians of BTD systems based on Laguerre 
	expansion}
	\label{alg2}
	\begin{algorithmic}[1]
		\Require Coefficient matrices $A_0,A_1,N_k (k = 1,2,\cdots,m),B,C$; the 
		parameters $\alpha$ and $p, q$.
		\Ensure The approximate Gramians $\tilde P$ and $\tilde Q$.
		\State Assemble the coefficient matrices $\Delta_i$ by \eqref{Delata_Coeff} 
		for $i=0, 1, \cdots, q-1$.
		\State Compute the factors $F_1$ and $G_1$ by \eqref{coeff_expp} \eqref{F1} 
		and \eqref{coeff_expQ} \eqref{G1}, respectively.
		\State Calculate the factors $F_j$ and $G_j$ by \eqref{coeff_pf} and 
		\eqref{Gi}, respectively, for $j=2, 3, \cdots, p$.
		\State Compute the low-rank approximate Gramians 
		\begin{equation*}
				\tilde P = F_1 F_1^{\top}+F_2 F_2^{\top}+\cdots+F_p F_p^{\top}, \,
				\tilde Q = G_1 G_1^{\top}+G_2 G_2^{\top}+\cdots+G_p G_p^{\top}. 
		\end{equation*}
		
	\end{algorithmic}
\end{algorithm}

\begin{remark}
	In Steps 2 and 3 of Algorithm \ref{alg2}, there is no need to assemble $\tilde 
	f_i$ and $\tilde g_i$ explicitly for $i=0, 1, \cdots, q-1$. The matrix-vector 
	products involved in $f_{j,i}$ and $g_{j,i}$ can be performed first, which 
	circumvents the matrix-matrix products contained in the definition of $\tilde 
	f_i$ and $\tilde g_i$, reducing the computational load dramatically in 
	large-scale settings. 
\end{remark}

\begin{remark}
    As pointed out in \cite{Benner2017add}, the 
	truncated Gramians, corresponding to a relatively small $p$, offer a reliable 
	approximation to Hankel singular values while are more efficient than full 
	Gramians for bilinear-like systems. In practice, the first two terms of the 
	summation in \eqref{eq8} is 
	commonly involved in the computation, and the value of $q\in [5, 10]$ typically 
	is sufficient  to provide an accurate approximation, which is also certified by 
	the simulation results provided in Section \ref{sec:sec-3.1}. In addition, the 
	uniform value of $q$ is not required for each term in the Laguerre expansion of 
	\eqref{eq8}, we adopt a consistent value in Algorithm \ref{alg2} for 
	simplicity.  
\end{remark}

\subsection{A BT procedure for BTD systems}

With the derived controllability and observability Gramians of BTD systems, a BT 
procedure can be designed to achieve the structure-preserving MOR. The conception 
of balanced transformation and balanced Gramians is presented. 
\begin{definition}
	Given BTD systems \eqref{eq1}, if there exists a nonsingular transformation 
	matrix $T$ such that the transformed systems, determined by 
	\{$\tilde{A}_0, \tilde{A}_1, \tilde{N}_1, \cdots,  \tilde{N}_m, \tilde{B}, 
	\tilde{C}\}$ with $\tilde{A}_0 = T^{-1} A_0 T,  \tilde{A}_1 = T^{-1} A_1 T , 
	\tilde{N}_k = 
	T^{-1} N_k T (k = 1,2,\cdots, m), \tilde{B}  = T^{-1} B, \tilde{C}  = C T$,
	possess Gramians $\tilde{P}=\tilde{Q}$ being diagonal, the matrix $T$ is called 
	a balanced transformation of \eqref{eq1}.
	The transformed systems $\{\tilde{A}_0, \tilde{A}_1, \tilde{N}_1, \cdots,  
	\tilde{N}_m, \tilde{B}, 
	\tilde{C}\}$ is referred as the balanced 
	system, and		
	$\tilde{P}$ and $\tilde{Q}$ are called balanced Gramians.
\end{definition}

Suppose that the system \eqref{eq1} has been balanced, and $\tilde{P} 
= \tilde{Q} = \operatorname{diag}\left\{\sigma_1, \sigma_2, \cdots, 
\sigma_n\right\}$ with $\sigma_i > \sigma_{i+1}$ for $i = 1, 2, \cdots, n-1$. We
partition the balanced system into the following form
$$
\tilde{A}_0=\left[\begin{array}{cc}
	A_0^{11} & A_0^{12} \\
	A_0^{21} & A_0^{22}
\end{array}\right],\tilde{A}_1=\left[\begin{array}{cc}
	A_1^{11} & A_1^{12} \\
	A_1^{21} & A_1^{22}
\end{array}\right],\tilde{N}_k=\left[\begin{array}{cc}
	N_k^{11} & N_k^{12} \\
	N_k^{21} & N_k^{22}
\end{array}\right](k = 1,2,\cdots, m), 
$$

$$
\tilde{B}=\left[\begin{array}{c}
	B^1 \\
	B^2
\end{array}\right], \quad \tilde{C}=\left[\begin{array}{ll}
	C^1 & C^2
\end{array}\right],
$$
where $A_0^{11},A_1^{11},N_k^{11} \in \mathbb{R}^{r \times r}\, \text{for}\, (k = 
1,2,\cdots, m), B^1 
\in \mathbb{R}^{r \times m}, C^1 \in \mathbb{R}^{l \times r} $. Here, $r$ is the 
reduced order, which can be determined by the approximate Hankel singular values 
with a given threshold. We truncate the states corresponding 
to the smaller singular values to obtain a ROM \eqref{eq2},
\ie, $\hat{A}_0 = A_0^{11}, \hat{A}_1 = A_1^{11}, \hat{N}_k = N_k^{11} (k = 
1,2,\cdots, m), \hat{B} = B^1, \hat{C} =  C^1$.

There are various schemes to produce a balance 
transformation. We employ the method based on Cholesky 
decomposition in order to avoid the direct calculation of the inverse matrix of $T$.
Specifically, the approximate Gramians $\tilde P$ and $\tilde Q$ of 
\eqref{eq1} are factorized first as 
\begin{equation*}\label{eqcholesky}
	\tilde P=U U^{\top}, \quad \tilde Q=L L^{\top},
\end{equation*}
then the SVD $L^{\top}  U=Z S Y \approx Z_r S_r Y_r^{\top}$ of the matrix $L^\top U 
$ is performed to construct the projection matrices $W_r = L Z_r S_r^{-1 / 2}$ and 
$V_r = U Y_r 
S_r^{-1 / 2}$. The resulting ROM \eqref{eq2} is obtained as follows 
\begin{equation*}
	\hat{A}_0 = W_r^{\top} A_0 V_r, \hat{A}_1 = W_r^{\top} A_1 V_r, \hat{N}_k = 
	W_r^{\top} N_k V_r (k = 1,2,\cdots, m), \hat{B} = W_r^{\top} B, \hat{C} =  C V_r.
\end{equation*}
The BT procedure for BTD systems is given in Algorithm 
\ref{alg1} based on approximate Gramians.

\begin{algorithm}[h]
	\caption{A BT procedure for BTD systems with approximate Gramians}
	\label{alg1}
	\begin{algorithmic}[1]
		\Require Coefficient matrices $A_0,A_1,N_k (k = 1,2,\cdots, m),B,C$.
		\Ensure Coefficient matrices $\hat{A}_0,\hat{A}_1,\hat{N}_k (k = 
		1,2,\cdots, m),\hat{B},\hat{C}$.
		\State Assemble the low-rank factors of approximate Gramians $\tilde P$ and 
		$\tilde Q$ by implementing Algorithm \ref{alg2}, which are defined as 
		\begin{equation*}
			F=[F_1, F_2, \cdots, F_p], \quad G=[G_1, G_2, \cdots, G_p]. 
		\end{equation*}
		\State Perform the SVD for the matrix $G^\top F$, that is $G^\top 
		F=ZSY^\top$. Determine a reduced order $r$ according to approximate Hankel 
		singular values, and assemble the matrices $Z_r = Z(:,1:r), S_r = 
		S(1:r,1:r), Y_r = (:,1:r)$.
		\State Construct the projection matrices $W_r=G Z_r S_r^{-1 / 2},V_r=F 
		Y_r 
		S_r^{-1 / 2}.$
		\State Calculate the coefficient matrices of ROMs
		$$
		\hat{A}_0 = W_r^\top A_0 V_r, \hat{A}_1 = W_r^\top A_1 V_r, 
		\hat{N}_k = 
		W_r^\top N_k V_r, \hat{B} = W_r^\top B, \hat{C} = C V_r.
		$$
	\end{algorithmic}
\end{algorithm}

\begin{remark}
	In Algorithm \ref{alg1}, we avoid to compute the approximate Gramians 
	explicitly, while extract the low-rank factors of Gramians directly from 
	Algorithm \ref{alg2}. This strategy not only circumvents the associated matrix 
	decomposition, such as Cholesky decomposition, involved in the standard BT 
	procedure, but also facilitates a lot the execution of SVD in Step 2 because of 
	the low-rank structure of the approximation. As a result, the presented BT 
	procedure is  
	efficient in large-scale settings. 
\end{remark}

\begin{remark}
	Although our discussion mainly focuses on BTD systems just with a single delay, 
	all results presented in the paper can be generalized to the case of multiple 
	delays with some proper modification. The main changes come from the fundamental 
	solution $K(t)$ in Definition \ref{def1} . For multiple delays $\tau_1, \cdots, 
	\tau_h$ with coefficient matrices $A_1, \cdots, A_h$, the main term associated 
	with the delays reads $\left(s I-A_0-\sum^h_{i=1} e^{-s \tau_i} A_i\right)^{-1}$ 
	for Gramians in the frequency domain, which can also be well approximated via 
	Laguerre expansion in a similar manner. Such an extension is somewhat trivial 
	and we omit the details in this paper for brevity. 
\end{remark}


\section{Numerical examples}\label{sec:sec-3.1}
We use two numerical examples to show the efficiency of the presented approach 
in 
this section. The simulation is conducted via Matlab (R2024b) on a laptop with 
Intel(R) Core(TM) Ultra 5 125H with 3.60 GHz and 16 GB RAM. For comparison, the 
moment-matching method based on the 
high-order Krylov subspace given in \citep{cheng2024high} is also carried 
out for BTD systems. The relative error of time response is adopted to access 
the accuracy of both methods. 
	
\begin{example}
		This example comes from a linear dynamical system with internal time 
		delay given in \cite{schulze2016data}. The system along with a bilinear 
		term is as follows
		\begin{equation*}
				\left\{\begin{array}{l}
				E \dot{x}(t)=A_0 x(t)+A_1 x(t-\tau)+N u(t) x(t)+B u(t), \\
				y(t)=C x(t), 
				\end{array}\right.
        \end{equation*}
		where the coefficient matrices is defined as 
		\begin{equation*}
				E=\theta I+T, \quad 
				A_0=\frac{1}{\tau}\left(\zeta^{-1}+1\right)(T-\theta I), \quad 
				A_1=\frac{\zeta^{-1}+1}{\zeta^{-1}-1} A_0. 
		\end{equation*}
		Specifically, $I$ is an identity matrix of order $n\times n$, and the 
		elements of $T$ are ones on the sub- and super-diagonal, in 
		the $(1, 1)$, and in the $(n, n)$ position and zeros everywhere else.  
		We use $N = 2J$ and $C = B^{\top}$, where $J$ is an $n$-dimensional 
		Jordan matrix with 
		zeros on 
		the main diagonal, and the first two elements of $B \in \mathbb{R}^n$ 
		are 1, while the rest are all zeros.

		In the simulation, we choose $\zeta = 0.01, \theta = 3, \tau = 3$, and 
		the dimension of the original system is $n = 500$. For our approach, the 
		parameter in Laguerre function is $\alpha=1$, and the first two terms of 
		the summation in \eqref{eq8} are used to get approximate Gramians. 
		Further, we adopt the first 5 terms in Laguerre polynomial expansion 
		\eqref{eq10} to simplify the involved numerical quadrature. The order of 
		resulting ROM is $r=6$. We consider the first two transfer functions to 
		conduct the moment-matching method, and the first 6 moments at the 
		origin are matching during the procedure of MOR, thereby leading to ROM 
		of order 6.

		\begin{figure}[htbp]
			\centering
			\includegraphics[width=1.00\linewidth]{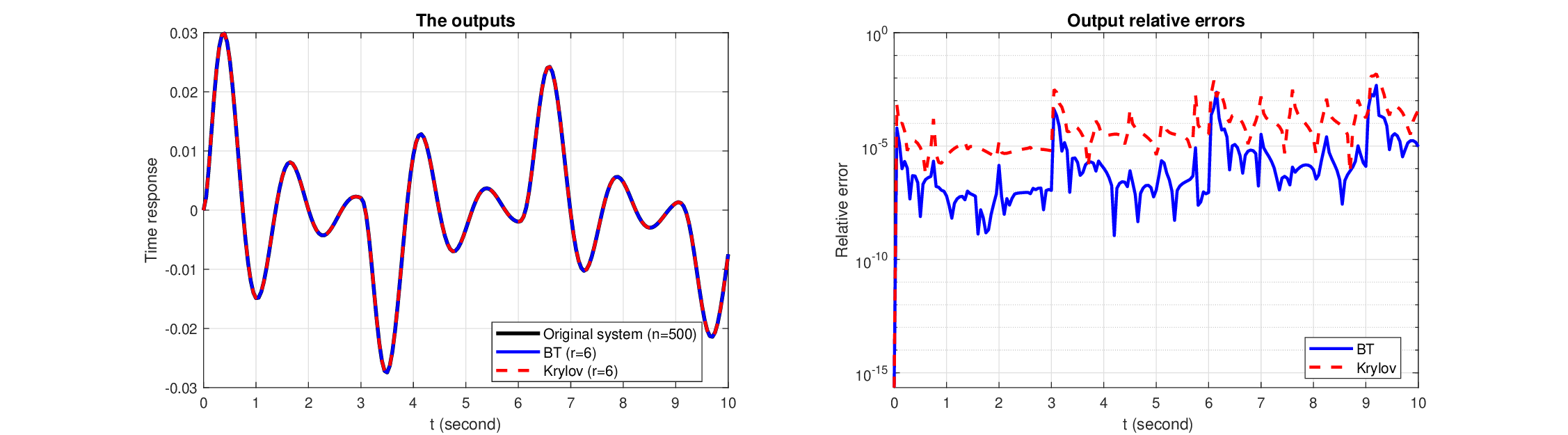}
			\caption {The dynamical responses and relative errors in the time 
			domain 
				in Example 1.}
			\label{fig1}
		\end{figure}
		
		\begin{figure}[htbp]
			\centering
			\includegraphics[width=0.45\linewidth]{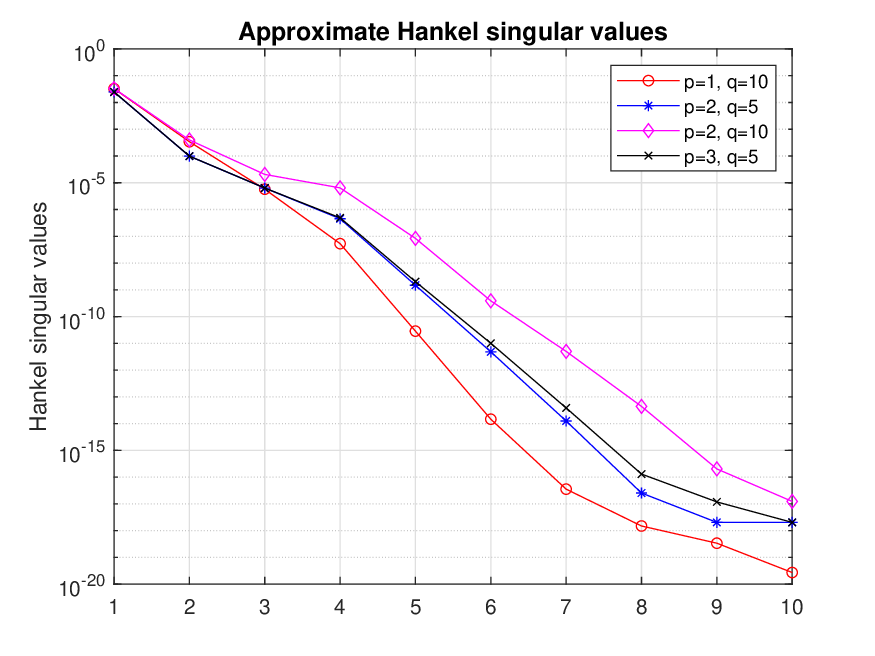}
			\caption {The approximate Hankel singular values for different 
				values in Example 1.}
			\label{fig2}
		\end{figure}
		
	\end{example}

		With the input function $u(t) = \mathrm{e}^{-t} \sin (5t)$, 
		Figure 
  \ref{fig1} 
  depicts the time response of ROMs and the associated relative errors. 
  The dynamical behavior of the original system is replicated faithfully 
  by ROMs for this example. The relative errors indicate that both ROMs 
  provide high-fidelity approximation to the original system, and the 
  proposed BT method possesses superior accuracy compared to 
  moment-matching method in the simulation. We also plot the first 10  
  approximate Hankel singular values of original systems with different 
  values in the approximation to Gramians. It is clear that the similarly 
  fast decay exhibits for each set of parameters. Note that taking the 
  first two 
  or three terms in \eqref{eq7} gives the comparable approximation to 
  Hankel singular values in this example, which admits the choice of 
  parameters for our approach in practice. 

	\begin{example}
		We consider a BTD system resulting from the variants of the following  
		partial differential equation with a delayed feedback 
		term
		\begin{equation*}
				\left\{\begin{array}{l}
				\frac{\partial u(x, t)}{\partial t}=\frac{\partial^2 u(x, 
				t)}{\partial x^2}-a_0(x) u(x, t)+a_1(x) u(\pi-x, t-1), \\
				u(0, t)=u(\pi, t)=0,
				\end{array}\right.
		\end{equation*}
 which describes a heating rod and frequently appears in 
 the field of MOR of time-delay systems \citep{peeters2013computing}, where 
 $a_0(x) = -2\sin(x), a_1(x) = 2\sin(x)$. The central difference with a 
 spatial step size of $h = \frac{\pi}{n+1}$ is adopted to perform the 
 discretisation in the interval $[0,\pi]$. As a result, $A_0$ is a 
 tridiagonal matrix, and $A_1$ is an anti-diagonal matrix.	The matrix $N$ 
 in the bilinear term is the identity matrix. The input matrix 
 $B^{\top} = [1\,\,1 \cdots 1] / \|[1\,\,1 \cdots 1]\|_2$. For
 the output matrix $C \in \mathbb{R}^{2 \times n}$, we have $C(1,:)=0.1$, 
 $C(2, 1)=C(2, n)=0.05$ and the other elements are $0$.

 \begin{figure}[htbp]
	  	\centering
	  	\includegraphics[width=0.98\linewidth]{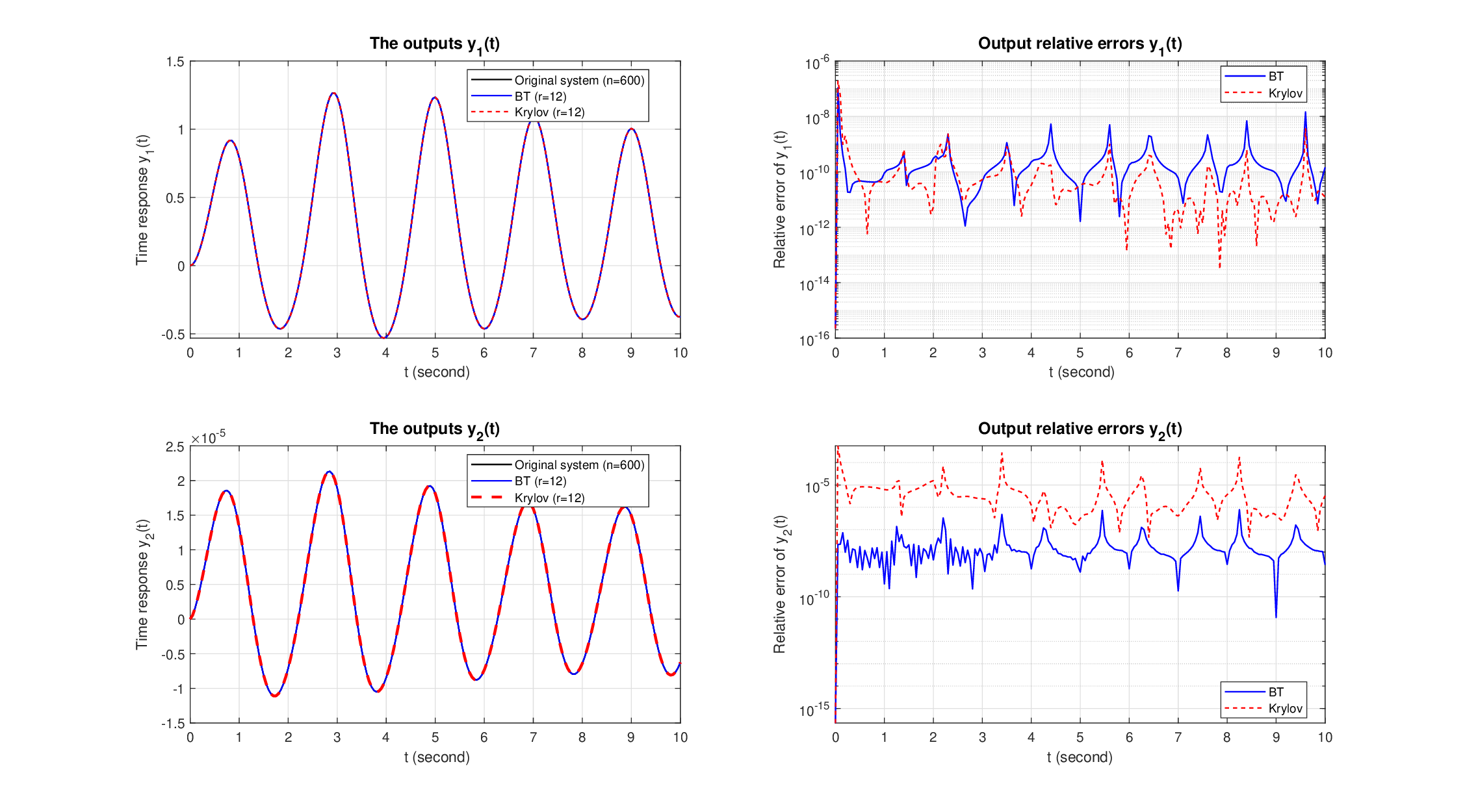}
	  	\caption {The dynamical responses and relative errors in the time domain 
	  		in Example 2.}
	  	\label{fig3}
 \end{figure}

 \begin{figure}[htbp]
	  	\centering
	  	\includegraphics[width=0.45\linewidth]{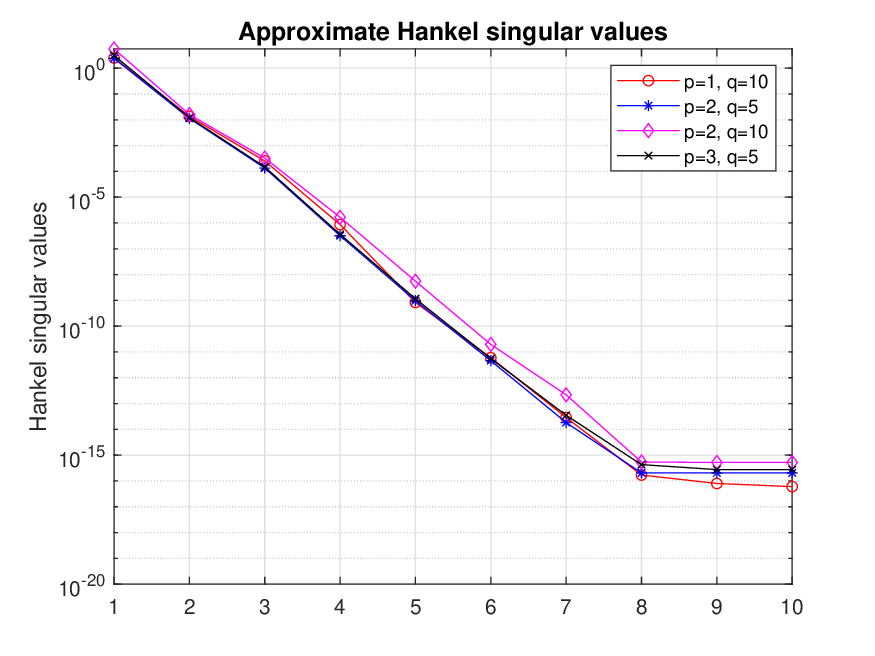}
	  	\caption {The approximate Hankel singular values for different 
	  		values in Example 2.}
	  	\label{fig4}
 \end{figure}
 
 We use $\tau = 2$ and $n = 600$ in the simulation. The first two terms of 
 \eqref{eq8} are involved in the approximation, and the first 5 terms are 
 involved in the truncated Laguerre expansion with the parameter 
 $\alpha=1$. The moment-matching method is also carried out to preserve the 
 first three moments for the first two transfer functions for this example. 
 The resulting ROMs are of order 12. Figure \ref{fig3} shows the 
 dynamical behaviors of the original system and ROMs, as well as the 
 relative errors when the systems are impulsed by the input function $u(t) 
 = \sin (\pi t)$. It is clear that both methods provide highly accurate 
 approximation to $y_1(t)$, while the proposed BT method performs much 
 better than the moment-matching method according to $y_2(t)$. In the 
 simulation, Figure \ref{fig4} shows that the different choice of 
 parameters in the numerical calculation of Gramians gives rise to almost 
 the same approximate Hankel singular values for this example, which 
 implies the accuracy of approximate Gramians. 				  
\end{example}

\section{Conclusions}\label{sec:sec-5}
We study a Gramian-based MOR method for BTD systems. The reachability 
and observability Gramians of BTD systems are given for the first time in the 
time domain and frequency domain based on the Volterra representation, which are 
compatible with the ones of bilinear systems and linear time-delay systems. A 
Laguerre expansion method is provided to calculate the approximate Gramians 
efficiently. A procedure of balanced truncation is proposed to produce 
structure-preserving ROMs. The simulation results demonstrate that the proposed 
BT procedure can provide high-precision approximation compared to the existing  
moment-matching method in some cases. There are still many issues that need to 
be further 
addressed. The stability and the error analysis of ROMs will be 
considered in the future work. 


\bibliography{reference}
\bibliographystyle{elsarticle-num}

\end{document}